\documentclass[12pt]{article}
\usepackage[T2A]{fontenc}
\usepackage{amsmath}
\usepackage[cp1251]{inputenc}
\begin{document}
\large

\begin{center}
{\bf IDENTIFICATION OF BOUNDARY CONDITIONS USING NATURAL
FREQUENCIES}
\end{center}

\vspace{0.2cm}

\centerline{AKHTYAMOV~A.~M.$^{a,}$\footnote{e-mail:
AkhtyamovAM@mail.ru}, MOUFTAKHOV~A.~V.$^{b,}$\footnote{e-mail:
muftahov@yahoo.com}}

\vspace{0.2cm}

$^a$ Department of Differential Equations,  Bashkir State
University, Ufa, Russia.

$^b$ Department of Mathematics and Statistics, Bar-Ilan University,
Ramat-Gan, Israel.

\vspace{0.2cm}

\date{}

\addtocounter{section}{1}
\begin{abstract}

The present investigation concerns a disc of varying thickness of
whose flexural stiffness $D$ varies with the radius $r$ according
to the law $D=D_0\, r^m$, where $D_0$ and $m$ are constants.
 The  problem of finding
 boundary conditions
for fastening this disc, which are inaccessible to direct
observation, from the natural frequencies of its axisymmetric
flexural oscillations is considered.
  The problem in question belongs to the class
of  inverse problems and is a completely natural problem of
identification of boundary conditions.
   The search for the unknown conditions
for fastening the disc
 is equivalent to finding
 the span of the vectors
of unknown   conditions coefficients. It is shown that this
inverse problem is well posed.

  Two theorems on the uniqueness and a theorem on stability
of the solution of this problem  are proved,  and a method for
establishing the unknown conditions for fastening the disc to the
walls is indicated. An approximate formula for determining the
 unknown conditions is obtained using
 first three natural frequencies.
 The method
of approximate  calculation of unknown boundary conditions is
explained  with the help of three examples of  different cases for
the fastening the disc (rigid clamping, free support, elastic
fixing).

\vspace{0.2cm}


\end{abstract}

\vspace{0.2cm}

{\bf 1. Introduction.} Discs are parts of turbines and other
various devices (see \cite{Campbell 24}--\cite{Vibrations 78}). If
it is impossible to observe the disc directly, the only source of
information about possible defects of its fastening can be the
natural frequencies of its flexural vibrations. The question
arises whether one would be able  to detect damage in disc
fastening by the natural frequencies of its symmetric flexural
vibrations. This paper gives and substantiates a positive answer
to this question.

 The problem in question
belongs to the class of inverse problems and is a completely
natural problem of identification of the boundary conditions.

Closely related formulations of the problem were proposed in
\cite{Kac 66,Qunli 90}. Contrary to this, in this paper it is not
the form of the domain or size of an object which are sought for
but the nature of fastening. The problem of determining a boundary
condition has been considered in \cite{Frikha 00}. However, as
data for finding the boundary conditions, we take not a set of
natural frequencies, but not condensation and inversion (as in
\cite{Frikha 00}).

Similarly formulated problems also occur in the spectral theory of
differential operators, where it is required to determine the
coefficients of a differential equation and the boundary
conditions using a set of eigenvalues (for more details, see
\cite{Borg 46}--\cite{Akhtyamov 99 DM}). However, as data for
finding the boundary conditions, we take  one spectrum  but not
several spectra or other additional spectral data (for example,
the spectral function, the Weyl function or the so-called
weighting numbers) that were used in \cite{Borg
46}--\cite{Akhtyamov 99 DM}. Moreover, the principal aim there was
to determine the coefficients in the equation and not in boundary
conditions. The aim of this paper is to determine  the boundary
conditions of the eigenvalue problem from its spectrum in the case
of a known differential equation.

The problem of determining a boundary condition using a finite set
of eigenvalues has been considered previously in
  \cite{Akhtyamov 99 DE}--\cite{Akhtyamov 01+}. In contrast to these
  papers, we think it is necessary to determine not
the  type of fastening of a string, membrane or  beam but the type
of disc fastening of varying
   thickness.

\vspace{0.2cm}

{\bf 2. Formulation of the direct problem.}  The problem of
axisymmetric flexural oscillations of disc of varying
   thickness
without a hole in the centre is reduced to following  differential
equation \cite{Timoshenko 40}
\begin{equation}\label{ipi 7 Conway 1} D\,\frac{\partial}{\partial r} \left(
\frac{\partial^2w}{\partial r^2}+\frac{1}{r}\, \frac{\partial
w}{\partial r}\right) + \frac{\partial D}{\partial r }\,\left(
\frac{\partial^2w}{\partial r^2}+\frac{\nu}{r}\, \frac{\partial
w}{\partial r}\right) =-\frac{1}{r}\, \int^r \rho\, h\,
\frac{\partial^2\, w}{\partial t^2}\, r\, dr,
\end{equation}
where $w$ is the deflection at radius $r$, $D$ is the flexural
rigitiy, $h$ is the thickness, $\rho$ is the density and $\nu$ is
Poisson's ratio.

If $D=D_0\, r^m$, where $D_0=E\, h_0^3/12(1-\nu^2)$, and m are
constants, the thickness $h$ is given by $h=h_0\, r^{m/3}$. Then
Equation (\ref{ipi 7 Conway 1})  becomes
$$
r^4\,\frac{\partial^4\, w}{\partial r^4}+r^3\, (2\, m+2)\,
\frac{\partial^3\, w}{\partial r^3}+r^2(m+\nu\, m-1+m^2)\,
\frac{\partial^2 w}{\partial r^2}
$$
$$
+r(1+\nu\, m^2-\nu\, m-m)\, \frac{\partial w}{\partial
r}=-\frac{12\, \rho\, (1-\nu^2)}{E\, h_0^2}\,
r^{(12-2m)/3}\,\frac{\partial^2 w}{\partial t^2}.
$$

For vibrations, we write
$$
w=u(r)\, e^{i\,\omega\, t}
$$
and hence obtain the following equation for $u(r)$:
\begin{equation}\label{ipi 7 Conway 4}
\begin{array}{c}
  r^4\,\dfrac{d^4\, u}{d r^4}+r^3\, (2\, m+2)\,
\dfrac{d^3\, u}{d r^3}+r^2(m+\nu\, m-1+m^2)\, \dfrac{d^2 u}{d r^2}
 \\[10pt]
  +r(1+\nu\, m^2-\nu\, m-m)\, \dfrac{d w}{dr}-\frac{12\, \rho\, (1-\nu^2)}{E\, h_0^2}\,
r^{(12-2m)/3}\, u=0.
\end{array}
\end{equation}

Although it is impossible to solve the problem for arbitrary $m$
exactly, some particular solutions in terms of Bessel functions
are found for a number of positive values of $m$ and corresponding
values of Poisson ratio $\nu$ in \cite{Conway 58}.

 For the sake of being definite,
we consider the particular case of Equation (\ref{ipi 7 Conway 4})
with $h=h_0\, r^{2/3}$, Poisson ratio $\nu= 1/9$, $\frac{320\,
\rho\, }{27\, E\, h_0^2}=1$,  and radius $R=1$.

 The problem of
axisymmetric flexural oscillations of a disc of varying
   thickness without a hole in the center for this particular case
 is reduced to following eigenvalue
problem \cite{Conway 58}
\begin{gather}
r^4\,\frac{d^{\, 4}\, u}{dr^4}+6\, r^3\,\frac{d^{\, 3}\,
u}{dr^3}+\frac{47\,r^2}{9}\,\frac{d^{\, 2}\, u}{dr^2}
-\frac{7\,r}{9}\, \frac{d\, u}{dr}-s^2\, r^{8/3}\, u= 0,
\label{ipi-eq1}
\\
U_1(u)=a_1\, \left( L_1\, u\right) _{r=1} + a_4\, \left( L_4\,
u\right) _{r=1} =0,\label{ipi-eq1a}
\\
U_2(u)=a_2\, \left( L_2\, u\right) _{r=1} + a_3\, \left( L_3\,
u\right) _{r=1}=0.\label{ipi-eq1b} \end{gather} Here
$s^2=\frac{12\, \rho\, \omega^2\, (1-\nu^2)}{E\, h_0^2}
=\frac{320\, \rho\, \omega^2}{27\, E\, h_0^2}=\omega^2$ is the
eigenvalue parameter, $E$ is the  elasticity modulus,  $\omega$ is
the natural frequencies parameter,
 $U_1(u),$ $U_2(u)$ are the  linear forms which characterize
conditions for fastening the plate to the walls (rigid clamping,
free support, free edge, floating fixing, elastic fixing),
\begin{gather*}
L_1\, u=u(r), \quad L_2\, u=\frac{du(r)}{dr}, \quad L_3\,
u=\frac{d^2 u(r)}{dr^2}+\frac{1}{9\, r}\,\frac{du}{dr},
\\
L_4\, u=\frac{d}{dr}\,\left[ \frac{d^2
u(r)}{dr^2}+\frac{1}{r}\,\frac{du}{dr}\right].
\end{gather*}

The solution to (\ref{ipi-eq1}) has form \cite{Conway 58}
\begin{gather*}\label{ipi-eq4}
  u=r^{-2/3}\, \left[ C_1\, I_{1}\left(\frac{3\, r^{2/3}\, {s}^{1/2}}{2}\right) +
   C_2\, J_{1}\left(\frac{3\, r^{2/3}\,  {s}^{1/2}}{2}\right) +\right.
  \\
\left.  +C_3\, Y_{1}\left(\frac{3\, r^{2/3}\,
{s}^{1/2}}{2}\right)
   +C_4\, Y_{1}\left(\frac{3\,i\, r^{2/3}\,  {s}^{1/2}}{2}\right)\right] .
\end{gather*}
Here $I_1(\cdot)$, $J_1(\cdot)$, $Y_1(\cdot)$ are the conventional
notations for first-order Bessel functions of real and imaginary
argument \cite{Watson 66}.

 For a continuous plate (without a hole in the center), the
constants $C_3=C_4=0.$

Instead of the  (\ref{ipi-eq1a}), (\ref{ipi-eq1b}), we can write
the conditions
\begin{equation}
 \ U_i(u)= \sum_{j=1}^4 a_{ij}\, \left( L_j\, u\right)
_{r=1},\qquad i=1,\, 2.\label{ipi-eq2}
\end{equation}
where $a_{11}=a_1,$ $a_{12}=0,$ $a_{13}=0,$ $a_{14}=a_4,$
$a_{21}=0,$ $a_{22}=a_2,$ $a_{23}=a_3,$ $a_{24}=0.$

We shall now formulate the direct eigenvalue problem
(\ref{ipi-eq1}), (\ref{ipi-eq2}): it is required to find the
unknown natural frequencies of the oscillations of the plate from
the
 linear forms $U_1(u),$ $U_2(u)$.

The natural frequencies $\omega_i$ are the corresponding positive
eigenvalues of problem (\ref{ipi-eq1}), (\ref{ipi-eq2}) (see
\cite{Vibrations 78}). The non-zero eigenvalues of problem
(\ref{ipi-eq1}), (\ref{ipi-eq2}) are the roots of the determinant
$$
   \Delta (s)=
\left|\begin{array}{cc}
  U_1(u_1) & U_1(u_2) \\
  U_2(u_1) & U_2(u_2)
\end{array}
\right| = M_{12}\, f_1(s)+M_{13}\, f_2(s)+M_{24}\, f_3(s)+M_{34}\,
f_4(s),
$$
 where\\ $u_1=u_1(r,s)=
r^{-2/3}\,I_1\left(\frac{3}{2}\, r^{2/3}\, {s}^{1/2}\right) , $
$u_2=u_2(r,s)= r^{-2/3}\,J_1\left(\frac{3}{2}\, r^{2/3}\,
{s}^{1/2}\right) $ are linearly independent solutions of Equation
(\ref{ipi-eq1});
\begin{equation}\label{ipi7 f}
\begin{array}{l}
f_1(s)=\left[\, L_1\, u_1(r,s)\cdot L_2\, u_2(r,s)-L_2\,
u_1(r,s)\cdot L_1\, u_2(r,s)\,\right]_{r=1} ,
\\
f_2(s)=\left[\, L_1\, u_1(r,s)\cdot L_3\, u_2(r,s)-L_3\,
u_1(r,s)\cdot L_1\, u_2(r,s)\,\right]_{r=1} ,
\\
f_3(s)=\left[\, L_4\, u_1(r,s)\cdot L_2\, u_2(r,s)-L_2\,
u_1(r,s)\cdot L_4\, u_2(r,s)\,\right]_{r=1} ,
\\
f_4(s)=\left[\, L_4\, u_1(r,s)\cdot L_3\, u_2(r,s)-L_3\,
u_1(r,s)\cdot L_4\, u_2(r,s)\,\right]_{r=1}.
\end{array}
\end{equation}

Thus, finding the three natural frequencies $\omega_i$ is
equivalent to founding of three roots $s_i$ of $\Delta (s)$.

For example, if boundary conditions (\ref{ipi-eq2})
 have the form
$$
U_1(u)=\left[ c\, L_1\, u(r)-L_4\, u(r)\right]_{r=1}=0, \quad
U_2(u)=\left[ L_2\, u(r)\right]_{r=1}=0
$$
(elastic clamping). For different $c$ we have the different
${s_i}$. See Table~1.
\begin{table}
$$
\begin{tabular}{|c||c|c|c|c|c|}
\hline
   & $c=10^{-5}$  & $c=1$ & $c=10$ & $c=10^2$ & $c=10^{100}$ \\
   \hline
  ${s_1}^{1/2}$ & 0.085457
 & 1.5178
 & 2.6517
 & 3.0663
 & 3.0739
 \\
${s_2}^{1/2}$ &  3.1122
 & 3.1145
 & 3.1561
 & 4.5575
 & 5.1995
  \\
  ${s_3}^{1/2}$ & 5.4634
 & 5.4651
 & 5.4813
 & 5.7412
 & 7.3054
 \\ \hline
\end{tabular}
$$
\caption{Dependence $c$ on $s_i.$}
\end{table}

Thus, knowing $c$ it is possible to find $s_i$ by standard
methods. The solution to this direct problem presents no
difficulties. The question arises whether one would be able to do
the reverse and find $c$ knowing $s_i$. In a broader sense it may
be stated as follows. Is it possible to derive unknown boundary
conditions with a knowledge of $s_i$? The answer to this question
is given in the next section.

\vspace{0.2cm}

{\bf 3. Formulation of the inverse problem.} The mathematical
(direct) problem is an eigenvalue problem for a homogeneous linear
fourth order equation, set up in the interval  $0<r<1$,
accompanied by two linear homogeneous boundary conditions at
$r=1$; the sought solutions must be bounded for $r=0$. The
boundary conditions depend on 4 scalar coefficients (namely 3
because of the homogeneity).

 Now we shall  formulate
the inverse of eigenvalue problem (\ref{ipi-eq1}),
(\ref{ipi-eq2}): it is required to find the unknown linear forms
$U_1(u),$ $U_2(u)$ from the natural frequencies of the disc
oscillations.

We shall denote the matrix, consisting of the coefficients
$a_{ij}$ of the forms $U_1(u)$ and  $U_2(u)$, by $A$ and its
minors by $M_{ij}$:
$$
A= \left\|\begin{array}{cccc} a_{11} & a_{12} & a_{13} & a_{14}
\\
a_{21} & a_{22} & a_{23} & a_{24}
\end{array}
\right\| , \qquad M_{ij}= \left|\begin{array}{cccc} a_{1i} &
a_{1j}
\\
a_{2i} & a_{2j}
\end{array}\right| .
$$

The search for the forms $U_1(u)$ and  $U_2(u)$ is equivalent to
finding the span $\left< {\bf a}_1,{\bf a}_2 \right>$ of the
vectors ${\bf a}_i=(a_{i1},a_{i2},a_{i3},a_{i4})^T$ $(i=1,2).$

 Hence, in terms of eigenvalue problem (\ref{ipi-eq1}),
(\ref{ipi-eq2}), the inverse problem constructed above should be
formulated as follows: the coefficients $a_{ij}$ of the forms
$U_1(u)$ and  $U_2(u)$ of problem (\ref{ipi-eq1}), (\ref{ipi-eq2})
are unknown, the rank of the matrix A to make up these
coefficients is equal to two, the non-zero natural frequencies
$s_k$ of problem (\ref{ipi-eq1}), (\ref{ipi-eq2}) are known and it
is required to find the span $\left< {\bf a}_1,{\bf a}_2 \right>$
of the vectors ${\bf a}_i=(a_{i1},a_{i2},a_{i3},a_{i4})^T$
$(i=1,2).$

\vspace{0.2cm}

{\bf 4. The uniqueness of the solution.}
 Together with problem (\ref{ipi-eq1}), (\ref{ipi-eq2}), let us
consider the following eigenvalue problem
\begin{gather}
r^4\,\frac{d^{\, 4}\, u}{dr^4}+6\, r^3\,\frac{d^{\, 3}\,
u}{dr^3}+\frac{47\,r^{\, 2}}{9}\,\frac{d^2\, u}{dr^2}
-\frac{7\,r}{9}\, \frac{d\, u}{dr}-s^2\, r^{8/3}\, u= 0,
\label{ipi-eq1v}
\\
\widetilde U_i(u)= \sum_{j=1}^4 b_{ij}\, \left( L_j\, u\right)
_{r=1},\qquad i=1,\, 2,\label{ipi-eq2v}
\end{gather}
where $b_{11}=b_1,$ $b_{12}=0,$ $b_{13}=0,$ $b_{14}=b_4,$
$b_{21}=0,$ $b_{22}=b_2,$ $b_{23}=b_3,$ $b_{24}=0.$

We denote the matrix composed of the coefficients $b_{ij}$ of the
forms $\widetilde U_1(u)$ and $\widetilde U_2(u)$ by $B$ and its
minors by $\widetilde M_{ij}$:
$$
B= \left\|\begin{array}{cccc} b_{11} & b_{12} & b_{13} & b_{14}
\\
b_{21} & b_{22} & b_{23} & b_{24}
\end{array}
\right\| , \qquad \widetilde M_{ij}= \left|\begin{array}{cccc}
b_{1i} & b_{1j}
\\
b_{2i} & b_{2j}
\end{array}\right| .
$$

The span of the vectors ${\bf
b}_i=(b_{i1},b_{i2},b_{i3},b_{i4})^T$ $(i=1,2)$ is denoted by
$\left< {\bf b}_1,{\bf b}_2\right>$.

{\it Theorem 1 (on the uniqueness of the solution of the inverse
problem).} Suppose the following conditions are satisfied
\begin{equation}\label{ipi 2 1}
  {\rm rank\, } A={\rm rank\, } B =2.
\end{equation}
If the non-zero eigenvalues $\{ s_k\}$ of problem (\ref{ipi-eq1}),
(\ref{ipi-eq2}) and the non-zero eigenvalues $\{ \widetilde s_k\}$
of problem (\ref{ipi-eq1v}), (\ref{ipi-eq2v}) are identical, with
account taken for their multiplicities, the spans $\left< {\bf
a}_1,{\bf a}_2\right>$ and $\left< {\bf b}_1,{\bf b}_2\right>$ are
also identical.

{\it Proof.}
 The non-zero eigenvalues of problem (\ref{ipi-eq1v}),
(\ref{ipi-eq2v}) are the roots of the determinant
$$
  \widetilde\Delta (s)= \left|\begin{array}{cc}
  \widetilde U_1(u_1) & \widetilde U_1(u_2) \\
  \widetilde U_2(u_1) & \widetilde U_2(u_2)
\end{array}
\right| = \widetilde M_{12}\, f_1(s)+\widetilde M_{13}\,
f_2(s)+\widetilde M_{24}\, f_3(s)+\widetilde M_{34}\, f_4(s). $$

In additional to the roots  identical to the non-zero eigenvalues
of the problems, the determinants $\Delta (s)$ and
$\widetilde\Delta (s)$ can also have the root $s=0$ of finite
multiplicity.

Since $\Delta (s)\not\equiv 0$, $\widetilde\Delta (s)\not\equiv 0$
are entire functions in s of order 1/2, it follows from Hadamard's
factorization theorem (see \cite{Gol'dberg 91}) that determinants
$\Delta (s)$  and $\widetilde\Delta (s)$ are connected by the
relation
$$
\Delta (s)\equiv C\, s^k \, \widetilde\Delta (s),
$$
where $k$ is a certain non-negative integer and $C$ is a certain
non-zero constant. From this we obtain the identity
\begin{equation}
\begin{array}{c}
(M_{12}-C\, s^k\,\widetilde M_{12}) \, f_1(s)+(M_{13}-C\,
s^k\,\widetilde M_{13})\, f_2(s)+\\
+(M_{24}-C\, s^k\,\widetilde M_{24})\, f_3(s)+(M_{34}-C\,
s^k\,\widetilde M_{34})\, f_4(s)\equiv 0. \end{array} \label{ipi
equiv}
\end{equation}
Note that the number $k$ in this identity is equal to zero.
Actually, let us assume the opposite: $k\ne 0$. Then the functions
$f_i(s)$ $(i=1,2,3,4)$ and also the same functions multiplied by
$s^k$ are linearly independent.

Indeed, using MAPLE, we get
\begin{eqnarray*}
f_1(s) &=&  - { \frac {27}{64}} \,s^{2} + { \frac {729}{131072}}
\,s^{4} - { \frac {6561}{671088640}} \,s^{6} +
{\frac {177147}{ 38482906972160}} \,s^{8} + {\rm O}(s^{10}),  \\
f_2(s) &=&  - {\displaystyle \frac {3}{16}} \,s^{2} +
{\displaystyle \frac {567}{32768}} \,s^{4} - {\displaystyle \frac
{9477}{167772160}} \,s^{6} + {\displaystyle \frac {373977}{
9620726743040}} \,s^{8}+ {\rm O}(s^{10}),\\
f_3(s) &=&\hspace{2cm} {\displaystyle \frac {81}{1024}} \,s^{4} -
{\displaystyle \frac {5103}{5242880}} \,s^{6} + {\displaystyle
\frac {216513}{150323855360}} \,s^{8}+ {\rm O}(s^{10}),\\
f_4(s) &=&\hspace{2cm} {\displaystyle \frac {9}{128}} \,s^{4} -
{\displaystyle \frac {81}{32768}} \,s^{6} + {\displaystyle \frac
{15309}{2684354560}} \,s^{8}+ {\rm O}(s^{10}).
\end{eqnarray*}
The determinant of the matrix $4\times 4$ of the coefficients at
$s^i$ ($i={2},\, 4,\, 6,\, 8$) in these representations of the
functions $f_j$ ($j=1,2,3,4$) is not equal $0$.  It now follows
that eight functions $f_j$, $s^k\, f_j$ ($j=1,2,\dots , 4$)  are
linearly independent for all $k\geq 8$.

Further, the determinant of the matrix $8\times 8$ of the
coefficients at $s^i$ ($i=4,\, 6,\, \dots \, , \, 18$) in the
series expansion of the functions $f_j$, $s^k\, f_j$ ($j=1,2,\dots
, 4$) is not equal $0$ under k=2, k=4 or k=6.  This means that
eight functions $f_j(s)$, $s^k\, f_j(s)$ ($j=1,2,\dots , 4$) are
linearly independent under k=2, k=4 or k=6.

This completes the the proof of linear independence of the
functions $f_j$, $s^k\, f_j$ ($j=1,2,\dots , 4$) for all $k\ne 0$.

From this and identity (\ref{ipi equiv}) we obtain
$$
M_{12}=M_{13}=M_{24}=M_{34}=\widetilde M_{12}=\widetilde
M_{13}=\widetilde M_{24}=\widetilde M_{34}=0,
$$
which, in combination with $M_{14}=M_{23}=0$, contradict condition
(\ref{ipi 2 1}) of the theorem.

Hence, $k=0$. From this and from identity (\ref{ipi equiv}), by
virtue of the linearly independence of the corresponding
functions, we obtain
$$
(M_{12},M_{13},M_{14},M_{23},M_{24},M_{34})^T= C\, (\widetilde
M_{12},\widetilde M_{13},\widetilde M_{14},\widetilde
M_{23},\widetilde M_{24},\widetilde M_{34})^T,
$$
which is equivalent to the proportionality of the bivectors ${\bf
a}_1\wedge {\bf a}_2$ and ${\bf b}_1\wedge {\bf b}_2$.

It is well-known \cite{Postnikov 79} that there is a natural
one-to-one correspondence between the classes of non-zero,
proportional bivectors and the two-dimensional subspaces of a
vector space. In this correspondence, a vector product ${\bf x}_1
\wedge {\bf x}_2$ of the vectors of its arbitrary basis ${\bf
x}_1$ and ${\bf x}_2$ corresponds to each subspace and a subspace
$\left< {\bf x}_1, {\bf x}_2\right>$ corresponds to each bivector
${\bf x}_1\wedge {\bf x}_2$. It therefore follows from the last
equation that $\left< {\bf a}_1, {\bf a}_2\right>$ which it was
required to prove.

\vspace{0.2cm}

{\bf 5. Exact solution.} It has been shown above that the problem
of finding the unknown linear forms $U_1(u)$ and $U_2(u)$ from the
natural frequencies of axisymmetric flexural oscillations of a
disc of varying thickness has a unique solution (in the sense that
the spans, composed of the coefficients of these linear forms, are
uniquely defined). The next question is how this solution can be
constructed.

This section deals with solving this problem and constructing
exact solution by the first three natural frequencies $\omega _i$.

Suppose  $s_1$,   $s_2$,  $s_3$ are the values corresponding to
the first three natural frequencies $\omega _i$. We substitute the
values $s_1$,   $s_2$,  $s_3$ into $\Delta(s)$ and obtain a system
of three homogeneous algebraic equations in the four unknows
$M_{12},M_{13},M_{24},M_{34}$
\begin{equation}\label{ipi7 sys}
  M_{12}\, f_1(s_i)+M_{13}\, f_2(s_i)+M_{24}\, f_4(s_i)+M_{34}\,
  f_4(s_i)=0,\qquad i=1,2,3.
\end{equation}

The resulting set of equations  has an infinite number of
solutions. If the resulting set has a rank of 3, the unknown
minors can be found in accurate to a coefficient.  The unknown
span and its basis are found from the minors using well-known
methods of algebraic geometry \cite{Hodge 94}. We also note
\cite{Akhtyamov 01} and \cite{Mouftakhov1}, in which detailed
method for solution of this problem was given.

The resulting set of equations has an infinite number of
solutions. It follows from the uniqueness theorem already been
proved that the unknown minors $M_{12},M_{13},M_{24},M_{34}$ can
be found accurate to a constant. Hence, the resulting system must
have a rank of 3 and the solution is determined accurate to a
constant multiplier.

If the minors $M_{12},M_{13},M_{24},M_{34}$ are found accurate to
a constant then the unknown span and its basis are found from the
minors using well-known methods of algebraic geometry \cite{Hodge
94}. We also note \cite{Akhtyamov 01} and \cite{Mouftakhov1}, in
which a detailed method for solution of this problem was given.

For example, if $M_{12}=1$, then
\begin{equation}\label{ipi pmm 36a1}
A= \left\|\begin{array}{cccc} 1 & 0 & 0 & -M_{24}
\\
0 & 1 & M_{13} & 0
\end{array}\right\| ;
\end{equation}
if $M_{13}=1$, then
\begin{equation}\label{ipi pmm 36a2}
A= \left\|\begin{array}{cccc} 1 & 0 & 0 & -M_{34}
\\
0 & M_{12} & 1 & 0
\end{array}\right\| .
\end{equation}


These reasons prove

{\it Theorem 2 (on the uniqueness of the solution of the inverse
problem).} {If the matrix of system (\ref{ipi7 sys}) has a rank of
3, the solution of the inverse problem of the reconstruction
boundary conditions~(\ref{ipi-eq1a}), (\ref{ipi-eq1b}) is unique.}

{\it Remark.} Theorem 2 is stronger than theorem 1. Theorem 2 use
only three natural frequencies for the reconstruction of boundary
conditions and
 not all natural frequencies as in
theorem 1.
  It follows from the uniqueness theorem 1 already
 proved that the unknown minors $M_{12},M_{13},M_{24},M_{34}$ can
be found accurate to a constant by all natural frequencies. To
prove that the unknown minors $M_{12},M_{13},M_{24},M_{34}$ can be
found accurate to a constant by three natural frequencies is very
easy problem in concrete cases. But in the common case to prove
that set (\ref{ipi7 sys}) has a rank of 3
 is very failed to do this. Yet our uniqueness theorem 1 suggests this result
 in the common case.

{\bf 6. Approximate solution.} Since small errors are possible
when measuring natural frequencies, the problem arises of finding
an algorithm for the approximate determination of the type of
fixing disc from  the first three natural frequencies  found with
a certain error.

If the values $s_i$ $(i=1,2,3)$ are approximately the same as the
first three exact eigenvalues of problem (\ref{ipi-eq1}),
(\ref{ipi-eq2}), and set (\ref{ipi7 sys})  has a rank equal to
three, the unknown minors can be determined accurate to a
coefficient.
 Further,  the problem arises of finding the
unknown span  from the values $  M_{12},  M_{13}, M_{24}, M_{34}$.
However, the values $  M_{12}$, $  M_{13}$, $  M_{24}$, $
  M_{34}$ cannot be minors of a matrix. So this problem
is not trivial. We must find minors $P_{12},P_{13},P_{24},P_{34}$
close to the values $  M_{12},  M_{13}, M_{24},  M_{34}$. This
problem is solved  with the help Lagrangian multiplier method and
algebraic geometry.

It is known from algebraic geometry that the numbers $P_{12}$,
$P_{13}$, $P_{14}$, $P_{23}$, $P_{24}$, $P_{34}$ are minors of
some matrix, if and only if the following condition is satisfied
$$ P_{12} P_{34} - P_{13} P_{24} + P_{14}P_{23} = 0.
$$ This condition is called Plucker condition (see \cite{Postnikov
79,Hodge 94}).

Using $P_{14} = P_{23} = 0$, we get
\begin{equation}\label{ipi pmm a2}
P_{12} P_{34} - P_{13} P_{24}  = 0.
\end{equation}

By definition, put $ \, x_1 = P_{12}, \, x_2 = P_{34},
 \, x_3 = P_{13},
 \, x_4 = -P_{24}.
$ Using this definition, we get Plucker condition
\begin{equation}\label{ipi pmm A3}
x_1\, x_2 + x_3\, x_4  = 0.
\end{equation}
It characterizes  a surface in the 4-dimensional space.

By definition, put $
 y_1 =   M_{12},$ $y_2 =   M_{34}$,
 $y_3 =   M_{13}$, $ y_4 = -  M_{24}$.

If Plucker condition for these numbers is realized, then $
  M_{12},$ $  M_{34}$,
 $  M_{13}$, $  M_{24}$ are minors of some matrix and
corresponding boundary conditions are found from the minors using
the methods in section 4.

If Plucker condition for numbers $   M_{12},$ $  M_{34}$,
 $  M_{13}$, $  M_{24}$ is not realized,  the
 required  minors
$P_{12}$, $P_{13}$,  $P_{24}$, $P_{34}$ are found with the help
Lagrangian multiplier method and algebraic geometry.

Indeed suppose $F(x_1,x_2,x_3,x_4,p)$ is Lagrange function
$$
(x_1-y_1)^2+(x_2-y_2)^2+(x_3-y_3)^2+(x_4-y_4)^2+2\, p\, (x_1\, x_2
+ x_3\, x_4),
$$
where $2\, p$ is the Lagrangian multiplier. If we find minimum of
the function $F(x_1,x_2,x_3,x_4,p )$, then we obtain the minors
$P_{12},P_{13},P_{24},P_{34}$ most close to the values $ M_{12},
M_{13},  M_{24},  M_{34}$.

The minimum of the function $F(x_1,x_2,x_3,x_4,p)$ is found from
equations
\begin{eqnarray}\label{ipi3 F'}
\left. \begin{array}{l} F_{x_1}'(x_1,x_2,x_3,x_4,p)=2\,
(x_1-y_1)+2\,p\, x_2=0,
\\
F_{x_2}'(x_1,x_2,x_3,x_4,p)=2\, (x_2-y_2)+2\, p\,x_1=0,
\\
F_{x_3}'(x_1,x_2,x_3,x_4,p)=2\, (x_3-y_3)+2\,p\, x_4=0,
\\
F_{x_4}'(x_1,x_2,x_3,x_4,p)=2\, (x_4-y_4)+2\,p\, x_3=0,
\end{array}\right\}
\\
F_{2\, p}'(x_1,x_2,x_3,x_4,p)=x_1\, x_2+x_3\, x_4=0. \label{ipi3
F'+}
\end{eqnarray}

By definition, put
$$
\begin{array}{c}
 X = (x_1, x_2 ,  x_3, x_4) ,\qquad Y = (y_1, y_2
, y_3, y_4),
\\
X^*=(x_2, x_1, x_4, x_3),\qquad Y^*=(y_2, y_1, y_4, y_3),
\\
(X, Y)=x_1\, y_1 + x_2\, y_2+x_3\, y_3 + x_4\, y_4. \end{array}
$$

With the preceding notations Equations (\ref{ipi3 F'}), (\ref{ipi3
F'+}) are identical to the following equations
\begin{gather}
Y = X + p\, X^*,\label{ipi pmm a4b}
\\
(X, X^*) = 0.\label{ipi pmm a4b+}
\end{gather}
  It follows from (\ref{ipi pmm a4b}) and (\ref{ipi pmm a4b+})
 that the vector $X^*$ is orthogonal to the vector $X$,
 and the vector  $ X = (x_1, x_2 ,  x_3, x_4) $ is orthogonal projection of
the vector  $Y = (y_1, y_2 , y_3, y_4) $ on the  surface (\ref{ipi
pmm A3}).

Having solved (\ref{ipi pmm a4b}) as set of linear equations with
the unknowns $x_1, x_2 ,  x_3, x_4$, we get
\begin{equation}\label{ipi pmm a4}
X = \frac1{1-p^2}\, (Y - p\, Y^*).
\end{equation}
From (\ref{ipi pmm a4}) it is easy to get
\begin{equation}\label{ipi pmm b4}
X^* = \frac1{1-p^2}\, (Y^* - p\, Y).
\end{equation}

Substituting (\ref{ipi pmm a4}) for $X$ and (\ref{ipi pmm b4}) for
$X^*$ in (\ref{ipi pmm a4b+}), we obtain
$$ (Y - p\, Y^*,\, Y^* - p\, Y) = 0.
$$



Notice that
$$
(Y,Y^*)\ne 0,\qquad (Y^*,Y^*)=(Y,Y), \qquad (Y^*,Y)=(Y,Y^*).
$$
Therefore,
$$
p^2 - 2\, p\, \frac {(Y,Y)}{(Y,Y^*)} + 1 = 0.
$$
This quadric equation has two roots
$$ p= \frac
{(Y,Y)\mp\sqrt{(Y,Y)^2 - (Y,Y^*)^2}}{(Y,Y^*)}.
$$
If $X$ is close to  $Y$, then $|p|<<1$ and thus we have
\begin{equation}\label{ipi pmm a5a}
p= \frac {(Y,Y)-\sqrt{(Y,Y)^2 - (Y,Y^*)^2}}{(Y,Y^*)}.
\end{equation}

The vector $X$ can be found using (\ref{ipi pmm a4}) and (\ref{ipi
pmm a5a}). The coordinates \linebreak
$P_{12},P_{13},P_{24},P_{34}$ of $X$ are minors of a matrix. This
matrix and corresponding boundary conditions are found from the
minors $P_{12},P_{13},P_{24},P_{34}$ using the methods in section
4.

\vspace{0.2cm}

{\bf 7.  Stability of the solution.}
 In this section we study continuity of the solution
 of the inverse problem with respect to $s_i$.
 It is shown that small perturbations of eigenvalues
 $s_i$ ($i=1,2,3$) lead to small perturbations of the
 boundary conditions.

Let $s_i$ ($i=1,2,3$) be eigenvalues of problem
(\ref{ipi-eq1}),(\ref{ipi-eq2}), \ and $\widetilde{s}_i$
($i=1,2,3$) are  such values that \vspace{0.2cm}
$$
|\widetilde{s}_i-s_i|<\delta<<1, \qquad i=1,2,3;$$ and $R$ is
such a number that
$$
|s_i|\leq R,\quad |\widetilde{s}_i|\leq R.
$$

\vspace{0.2cm}

{\it Theorem 3.} Suppose that one of the third-order minors of the
matrix $\|f_k(s_i)\|_{3\times 4}$ is substantially non-zero. If $
|\widetilde{s}_i-s_i|<\delta<<1,$ then the boundary conditions
 of problem (\ref{ipi-eq1v}),(\ref{ipi-eq2v}) are close to
 the boundary conditions
 of problem (\ref{ipi-eq1}),(\ref{ipi-eq2}).

\vspace{0.2cm}

{\it Proof.}
 By virtue of continuity
of the functions (\ref{ipi7 f}), we have
\begin{equation}\label{ip2 nepr f}
|f_k(\widetilde{s}_i)-f_k(s_i)|<M\,\delta, \qquad i=1,2,3,\quad
k={1,2,3,4},
\end{equation}
where $M=\max\limits_{k=1,2,3,4}M_k,$ \ $M_k=\max\limits_{|z|<R}
f_k(z)$.

First, let us prove that
$$
|\widetilde{M}_{ij}-{M}_{ij}|<C_1\, \delta .
$$

Recall that one of the third-order minors of the matrix
$\|f_k(s_i)\|_{3\times 4}$ is substantially non-zero by the
hypotheses of Theorem 3, and so
\begin{equation}\label{ip2 f_k}
 {\rm rank}\,
\|f_k(s_i)\|_{3\times 4}=3, \qquad i={1,2,3},\quad k=1,2,3,4.
\end{equation}

To be precise, assume that
\begin{equation}\label{ip2 1f_k}
 \mod(\det
\|f_k(s_i)\|_{3\times 3})>>0,\qquad i={1,2,3},\quad k=1,2,3.
\end{equation}
It follows from (\ref{ip2 nepr f}) that
\begin{equation}\label{ip2 2f_k}
 \mod(\det
\|f_k(\widetilde{s}_i)\|_{3\times 3})>>0,\qquad i={1,2,3},\quad
k=1,2,3.
\end{equation}

It can be shown by direct calculations that  set  (\ref{ipi7 sys})
has the solution
\begin{equation}\label{ip2 2M}
\begin{array}{l}
  M_{12}=\det(\|f_k(s_i)\|_{3\times 3}), \qquad i={1,2,3},\quad k=2,3,4, \\
    M_{13}=\det(\|f_k(s_i)\|_{3\times 3}), \qquad i={1,2,3},\quad
    k=1,3,4,
 \\
      M_{24}=\det(\|f_k(s_i)\|_{3\times 3}), \qquad i={1,2,3},\quad k=1,2,4, \\
      M_{34}=\det(\|f_k(s_i)\|_{3\times 3}), \qquad i={1,2,3},\quad k=1,2,3.
\end{array}
\end{equation}
 Like previously, the set
 \begin{equation}\label{IP ipi pmm 36v}
\widetilde{M}_{12}\,f_1(\widetilde{s}_i ) + \widetilde{M}_{13} \,
f_2(\widetilde{s}_i ) + \widetilde{M}_{24} \, f_3(\widetilde{s}_i
) +\widetilde{M}_{34} \, f_4(\widetilde{s}_i )=0, \qquad i=1,2,3.
\end{equation}
has  solution
\begin{equation}\label{ip2 2Mv}
\begin{array}{l}
  \widetilde{M}_{12}=\det(\|f_k(\widetilde{s}_i)\|_{3\times 3}), \qquad i={1,2,3},\quad k=2,3,4, \\
    \widetilde{M}_{13}=\det(\|f_k(\widetilde{s}_i)\|_{3\times 3}), \qquad i={1,2,3},\quad
    k=1,3,4,
 \\
      \widetilde{M}_{24}=\det(\|f_k(\widetilde{s}_i)\|_{3\times 3}), \qquad i={1,2,3},\quad k=1,2,4, \\
      \widetilde{M}_{34}=\det(\|f_k(\widetilde{s}_i)\|_{3\times 3}), \qquad i={1,2,3},\quad k=1,2,3.
\end{array}
\end{equation}

By definition, put
\begin{gather*}
\Delta_{ik}=f_k(\widetilde{s}_i)-f_k(s_i),\quad
\Delta_{k}=\|\Delta_{ik}\|_{3\times 1},\quad
F_{k}=\|f_k(s_i)\|_{3\times 1}, \\ ( i=1,2,3,\quad k={1,2,3,4}).
\end{gather*}

In the new notations, we have
\begin{eqnarray*}
&&|\widetilde{M}_{12}-{M}_{12}|=
\mod\big(\det(\|f_k(\widetilde{s}_i)\|_{3\times
3})-\det(\|f_k({s}_i)\|_{3\times 3})\big)=
\\
 &&= \mod\big(\det\|F_1+\Delta_{1},\, F_2+\Delta_{2},\, F_3+\Delta_{3}\|-
\det\|F_1,\, F_2,\, F_3\|\big)=
 \\
  &&= \mod\big(\det\|F_1,\, F_2+\Delta_{2},\,
  F_3+\Delta_{3}\|+
\\
&&\hspace{1cm}+\det\|\Delta_{1},\, F_2+\Delta_{2},\,
  F_3+\Delta_{3}\|-
 \\
 &&\hspace{1cm} -\det\|F_1,\, F_2,\, F_3\|\big)=
\\
  &&= \mod\big(\det\|F_1,\, F_2,\,  F_3\|+
\det\|\Delta_{1},\, F_2,\, F_3\|+
  \\
   &&\hspace{1cm} +\dots +\det\|\Delta_{1},\, \Delta_{2},\,\Delta_{3}\|-
\det\|F_1,\, F_2,\,  F_3\|\big)=
\\
 && =\mod\big(\det\|\Delta_{1},\, F_2,\,  F_3\|+
  \dots +\det\|\Delta_{1},\, \Delta_{2},\,\Delta_{3}\|\big).
 \end{eqnarray*}

It follows from (\ref{ip2 nepr f}) that
\begin{equation}\label{ip2 2ner}
 |\Delta_{ik}|<M\,\delta.
\end{equation}
Calculating the determinant and using (\ref{ip2 2ner}), we get
$$\mod\big(\det\|\Delta_{1},\, F_2,\,
F_3\|\big)<6\cdot M^3\cdot \delta.
  $$
Similarly, \ \ $
\mod\big(\det\|\Delta_{1},\,\Delta_{2},\,F_3\|\big)<6\cdot
M^3\cdot \delta^2\leq 6\cdot M^3\cdot \delta.
  $

Arguing as above, we see that $$ \mod\big(\det\|\Delta_{1},\,
F_2,\, F_3\|+
  \dots +\det\|\Delta_{1},\, \Delta_{2},\,\Delta_{3}\|\big)<
 7\cdot 6\cdot M^3\cdot \delta .
  $$
By definition, put $C_1=7\cdot 6\cdot M^3$. Then,
$$
|\widetilde{M}_{12}-{M}_{12}|=\mod\big(\det(\|f_k(\widetilde{s}_i)\|_{3\times
3})-\det(\|f_k({s}_i)\|_{3\times 3})\big)<C_1\,\delta .
$$

Continuing this line of reasoning, we see that
$$
|\widetilde{M}_{ij}-{M}_{ij}|<C_1\,\delta .
$$

This implies that
\begin{equation}\label{ip2 2Y-Y}
 |\widetilde{Y}-Y|
=\sqrt{\sum_{k=1}^4 (\widetilde{y}_k-y_k)^2}<2\, C_1\,\delta ,
\end{equation}
where \begin{gather*}
 Y=(M_{13},\, -M_{24},\, M_{14},\,
M_{23})^T, \qquad \widetilde{Y}=(\widetilde{M}_{13},\,
-\widetilde{M}_{24},\, \widetilde{M}_{14},\,
\widetilde{M}_{23})^T.
\end{gather*}
 By (\ref{ip2 1f_k}), (\ref{ip2 2f_k}) this means that
$\widetilde{Y}$ is close to $Y$.

 Now let us prove that
${\widetilde{X}=(\widetilde{P}_{12},\, \widetilde{P}_{34},\,
\widetilde{P}_{13},\, -\widetilde{P}_{24})^T}$ close to \linebreak
$X=(P_{12},\, P_{34},\, P_{13},\, -P_{24})^T$.

Combining  (\ref{ipi pmm a5a}) and (\ref{ip2 2Y-Y}), we see that
$$
\widetilde{p}= \frac
{(\widetilde{Y},\widetilde{Y})-\sqrt{(\widetilde{Y},\widetilde{Y})^2
-
(\widetilde{Y},\widetilde{Y}^*)^2}}{(\widetilde{Y},\widetilde{Y}^*)}
$$
close to $p$.
 This means that
 \ $
 \widetilde{X}= (\widetilde{Y} - \widetilde{p}\,
 \widetilde{Y}^*)/(1-\widetilde{p}^2)
 $ \
close to (\ref{ipi pmm a4}).

 Finally, let us prove that the boundary conditions
 of problem (\ref{ipi-eq1v}),(\ref{ipi-eq2v}) are close to
 the boundary conditions
 of problem
  (\ref{ipi-eq1}),(\ref{ipi-eq2}),
 where
 $$
 |\widetilde{s}_i-s_i|<\delta<<1.
$$

We know already (see \cite{Mouftakhov1}) that the coefficients of
the boundary conditions are the minors of matrix $A$ as in
(\ref{ipi pmm 36a1}).
 This means that if $\widetilde{X}$ is close
to $X$, the boundary conditions
 of problem (\ref{ipi-eq1v}),(\ref{ipi-eq2v}) are close to
 the boundary conditions
 of problem (\ref{ipi-eq1}),(\ref{ipi-eq2}).

This completes the proof.

\vspace{0.2cm}

Computer calculations confirm the stability of the solution
 of the inverse problem. The order of error often
hardly different from the error in the closeness of values of
$\widetilde\lambda_i$ and $\lambda_i$ and only in some cases it
can be deteriorated by four orders of magnitude. So the
measurement accuracy of instruments to measure natural frequencies
must exceed accuracy to measure boundary conditions by four orders
of magnitude.

It follows from theorem~2 and 3  that the inverse problem is well
posed, since its solution exists, is unique and
 continuous with respect to $s_i$ $(i=1,2,3)$.

\vspace{0.2cm}

{\bf 8. Examples.}  We use dimensionless variables in the
numerical examples.

{\it Example 1.} If $\widetilde{ s_1}^{\; 1/2}=3.0739$, $
\widetilde{s_2}^{\; 1/2}=5.1995$, $\widetilde{  s_3}^{\;
1/2}=7.3054$ correspond to the first three natural frequencies
$\omega_i$
 determined using instruments for measuring the natural frequencies
 with an accuracy of $10^{-4}$, then the solution of set
 (\ref{ipi7 sys}), accurate to a constant, has the form
$$
M_{12}=645330\, C,\quad M_{13} = 19.947\, C, \quad M_{24} =
2.4157\, C,\quad M_{34}=C.
$$
Using (\ref{ipi pmm a4}) and (\ref{ipi pmm a5a}), we get
$$\begin{array}{cccccc}
P_{12} & = & 645330\, C,\quad             &  P_{13} & = &
19.94703753\, C,\quad
\\
P_{24} & = & 2.416009841\, C,\quad & P_{34} & = & 0.000749141\, C.
\quad
\end{array}
$$

Suppose  $C=1/{P_{12}}$; then from (\ref{ipi pmm 36a1}), we obtain
$$
A= \left\|\begin{array}{cccc} 1 & 0 & 0 & -3.75565\cdot 10^{-5}
\\
0 & 1 & 0.00031  & 0
\end{array}\right\|.
$$

Note that the  numbers $\widetilde{  s_1}^{\; 1/2}=3.0739$,
$\widetilde{ s_2}^{\; 1/2}=5.1995$, $\widetilde{  s_3}^{\;
1/2}=7.3054$ presented above are almost the same as the first
three exact values ${s_1}^{1/2}$, ${s_2}^{1/2}$, ${s_3}^{1/2}$
corresponding to rigid clamping. This means that the unknown disc
fastening inaccessible to direct observation has been correctly
determined.

\vspace{0.2cm}

 {\it Example 2.} If $\widetilde{  s_1}^{\; 1/2}=1.8312$,
$  \widetilde{s_2}^{\; 1/2}=4.4629$, $\widetilde{  s_3}^{\;
1/2}=6.6502$ correspond to the first three natural frequencies
$\omega_i$
 determined by a frequency meter
 with an accuracy of $10^{-4}$, then the solution of set
 (\ref{ipi7 sys}),  accurate to a constant, has the form
$$
M_{12}=-9.31\, C,\; M_{13} = 201420\, C, \; M_{24} = -2.6702\,
C,\; M_{34}=C.
$$
Using (\ref{ipi pmm a4}) and (\ref{ipi pmm a5a}), we get
$$\begin{array}{cccccc}
P_{12} & = & -9.310013258\, C,\quad             &  P_{13} & = &
201420\, C,\quad
\\
P_{24} & = & -4.62276\cdot 10^{-5}\, C,\quad & P_{34} & = &
1.00012342\, C. \quad
\end{array}
$$


Suppose  $C=1/{P_{13}}$; then from  (\ref{ipi pmm 36a2}), we
obtain
$$
A= \left\|\begin{array}{cccc}
 1 & 0 & 0 & -4.96536\cdot 10^{-5}
\\
0 &  -4.62219\cdot 10^{-5} & 1 & 0
\end{array}\right\|.
$$

Note that the  numbers $\widetilde{ s_1}^{\; 1/2}$,
$\widetilde{s_2}^{\; 1/2}$, $\widetilde{s_3}^{\; 1/2}$ presented
above are almost the same as the first three exact values
${s_1}^{1/2}$, ${s_2}^{1/2}$, ${s_3}^{1/2}$ corresponding to the
free support. This means that the unknown  disc fastening
inaccessible to direct observation has been correctly determined.

\vspace{0.2cm}

 {\it Example 3.}
  If $\widetilde{s_1}^{\; 1/2}=1.5178$,
$  \widetilde{s_2}^{\; 1/2}=3.1145$, $\widetilde{  s_3}^{\;
1/2}=5.4651$ correspond to the first three natural frequencies
$\omega_i$
 determined by means of instruments for measuring  natural frequencies
 with an accuracy of $10^{-4}$, then the solution of set
 (\ref{ipi7 sys}),  accurate to a constant, has the form
$$
M_{12}=-140260 \, C,\quad M_{13} = 154.74 \, C, \quad M_{24} =
140150 \, C,\quad M_{34}=C.
$$

Using (\ref{ipi pmm a4}) and (\ref{ipi pmm a5a}), we get
$$\begin{array}{cccccc}
P_{12} & = & -140260 \, C,\quad             &  P_{13} & = & 76.931
\,
C,\quad \\
P_{24} & = & 140150 \, C,\quad & P_{34} & = & -76.870 \, C. \quad
\end{array}
$$


Suppose  $C=1/{P_{12}}$; then from (\ref{ipi pmm 36a1}), we obtain
$$
A= \left\|\begin{array}{cccc} 1 & 0 & 0 & -0.99922

\\
0 & 1 & -0.00054849
  & 0
\end{array}\right\|.
$$

Note that the  numbers $\widetilde{s_1}^{\; 1/2}$,
$\widetilde{s_2}^{\; 1/2}$, $\widetilde{s_3}^{\; 1/2}$ presented
above are almost the same as the first three exact values
${s_1}^{1/2}$, ${s_2}^{1/2}$, ${s_3}^{1/2}$, which correspond to
elastic fixing
with matrix $$ A=\left\| \begin{array}{cccc} 1&0&0&-1 \\
0&1&0&0 \end{array} \right\| .$$
 This means that the unknown
 disc fastening inaccessible to direct observation has been satisfactorily
determined.

\vspace{0.2cm}

Thus, the form of the disc fastening of varying thickness can be
determined from the first natural frequencies measured by special
instruments.

Note that we choose such a particular variation law for the
thickness of the disc to be precise. If the variation law for the
thickness of the disc is different from that adopted in this
study, then the mathematical formulation and the proposed
procedures of solving of the direct and inverse problems remain
valid. In this case we must substitute linearly independent
solutions  of corresponding differential equation for $u_1$ and
$u_2$ in (\ref{ipi7 f}), whose are not singular at $r=0$,
 and the corresponding linearly forms  for $L_1$, $L_2$, $L_3$,
$L_4$  in (\ref{ipi-eq2}).

Direct problems on  hydroelasticity and aeroelasticity  are
considered in \cite{Ilgamov 88,Ilgamov 91}. Similar inverse
problems can be solved by means of the method proposed in this
paper.

\newpage

{\bf 9. Acknowledgements.}

The authors are grateful to professor M.~A.~Ilgamov  for useful
discussions of engineering aspects.

This research was partially supported by the Russian Foundation
for Basic Research (01-01-00996), Ministry of Education of Russia
(E02-1.0-77), Emmy Noether Research Institute for Mathematics, the
Minerva Foundation of Germany, the Excellency Center "Group
Theoretic Methods in the Study of Algebraic Varieties"
 of the Israel Science Foundation, and by EAGER
(European Network in Algebraic Geometry).


\normalsize

\begin{thebibliography}{99}


\bibitem{Campbell 24}
 W. Campbell, The protection of steam-turbine disk wheels from axial
vibration, {\it Transactions of ASME}, {\bf 46}, 31--160 (1924).




\vspace{-0.3cm}

\bibitem{Levin 37}
 A.V. Levin, Vibration of disks, {\it Journal of technical physics},
{\bf 7}, No.~17, 1739--1753 (1937). \vspace{-0.3cm}

\bibitem{Strutt 29}
 W. Strutt (Lord Rayleigh), {\it The theory of sound.} 2d ed. Dover Publications, New
 York, N.Y.,
 1945,  V.~1, p.~xlii+480.

\vspace{-0.3cm}

\bibitem{Timoshenko 55}
{S. Timoshenko}, {\it Vibration Problems in Engineering},  D. Van
Nostrand Company, New York, 1937, p.~ix+470.

\vspace{-0.3cm}

\bibitem{Timoshenko 40}
{S. Timoshenko}, {\it Theory of Plates and Shells}, D. Van
Nostrand Company, New York, 1940, p.~283.

\vspace{-0.3cm}

\bibitem{Vibrations 78}
V.~V.~Bolotin (Ed.), {\it Vibrations in Engineering: A Handbok,
Vol.~1, Oscillations of Linear Systems},  Mashinostroenie, Moscow,
1978, p.~352.

\vspace{-0.3cm}

\bibitem{Kac 66} M.~Kac, Can one hear the shape
of a drum?,
 {\it Amer. Math. Monthly},  {\bf 73}, No. 4, 1--23 (1966).

\vspace{-0.3cm} \bibitem{Qunli 90} W. U. Qunli,  F. Fricke,
Determination of the size of an object and its location in a
cavity by eigenfrequency shifts, {\it Nat. Conf. Publ./ Inst. Eng.
Austral}, No. 9, 329--333,  (1990).

\vspace{-0.3cm} \bibitem{Frikha 00} {S. Frikha, G. Coffignal,
J.~L. Trolle}, Boundary condition identification using
condensation and inversion, {\it J. Sound and Vib.}, {\bf 233},
No.~3, 495--514 (2000) .

 \vspace{-0.3cm}


\bibitem{Borg 46}  G. Borg,
 {Eine umkehrung der Sturm---Liouvilleschen
  eigenwertanfgabe. Bestimmung der Differentialgleichung
  durch die Eigenwarte}, {\it  Acta Math.}, {\bf 78}, No.~1, 1--96  (1946).

\vspace{-0.3cm}

\bibitem{Marchenko 77} V.~A.~Marchenko, {\it Sturm-Liouville Operators and
their Applications}, Naukova Dumka, Kiev, 1977, p.~331; English
transl.: Birkh\"auser, Basel, 1986, p.~xii+367.

\vspace{-0.3cm}

\bibitem{Levitan 87} B. M. Levitan, {\it Inverse Sturm-Liouville Problems},
Nauka, Moscow, 1984, p.~240; English transl., VNU Science Press,
Zeist, 1987, p.~x+240.

\vspace{-0.3cm}

\bibitem{Poshel 87}
J.~P\"oshel and E.~Trubowitz, {\it Inverse Spectral Theory},
Academic Press, Boston, MA, 1987, p.~x+192.

\vspace{-0.3cm}

\bibitem{Yurko 00}
 V. A. Yurko, {\it Inverse Spectral Problems for Linear Differential Operators
and their Applications}, Gordon and Breach, New York, 2000. p.
253.

\vspace{-0.3cm}

\bibitem{Akhtyamov 99 DM}
V.~A. Sadovnichii, Ya.~T. Sultanaev, A.~M. Akhtyamov, Analogues of
Borg's uniqueness theorem in the case of
 nonseparated boundary conditions, {\it
 Doklady Mathematics},  {\bf 60}, No.~1, 115--117 (1999).


\vspace{-0.3cm}

\bibitem{Akhtyamov 99 DE}
A.~M. Akhtyamov,  Determination of the boundary condition on the
basis of a finite set of eigenvalues, {\it Differential
equations}, {\bf 35}, Part~8, 1141--1143 (1999).


\vspace{-0.3cm}

\bibitem{Akhtyamov 01}
I. Sh. Akhatov,  A. M. Akhtyamov,   Determination of the form of
attachment of the rod using the natural frequencies of its
flexural oscillations, {\it J. Appl. Maths Mechs}, {\bf 65},
No.~2, 283--290 (2001),

\vspace{-0.3cm}

\bibitem{Akhtyamov 01+}
 Akhtyamov~A.~M., Recognition of fastening the annular membrane
 from the natural frequencies of its oscillations,
 {\it Transactions of RANS, seies MMMIC},  {\bf 5}, No.~3, 103--110 (2001).



\vspace{-0.3cm}

\bibitem{Conway 58}  H. D. Conway, Some special solutions for
flexurel vibrations of discs of varying thickness, {\it Ing.
Arch.} {\bf 26}, No. 6, 408--410 (1958).

\vspace{-0.3cm}


\bibitem{Watson 66} G.~N. Watson, {\it A Treatise on the Theory of
Bessel Functions}, Cambridge Univ. Press, Cambridge, UK, 1995,
p.~vii+804.

\vspace{-0.3cm}

\bibitem{Gol'dberg 91}
B.~Ya.~Levin, {\it Distribution of zeros of entire functions},
Gostekhizdat,  Moscow, 1956. p. 632; English transl.: Amer. Math.
Soc., Providence, R.~I., 1980, p.~524.

\vspace{-0.3cm}

\bibitem{Postnikov 79} M. M. Postnikov, {\it Linear Algebra and
Differential Geometry},  Nauka, Moscow, 1979, p.~312; English
transl.: Moscow, MIR, 1982, p. 319.

\vspace{-0.3cm}

\bibitem{Hodge 94}
W.~V.~D.~Hodge, D.~Pedoe, {\it Methods of Algebraic Geometry},
Cambridge Univ. Press., Cambridge, UK, 1994, p.~viii+440.

\vspace{-0.3cm}

\bibitem{Mouftakhov1} A. V. Mouftakhov,
 On the reconstruction of the matrix from its minors, {\it 31-Mar-2003, MPS: Pure mathematics/03040001} (2003).

\vspace{-0.3cm}

\bibitem{Ilgamov 88} Dowell E.~H., Ilgamov M.~A. {\it Studies in
Nonlinear Aeroelasticity}, Springer Verlag, New York - Tokyo,
1988, p.~456.

\vspace{-0.3cm}

\bibitem{Ilgamov 91} Ilgamov~M.~A. {\it Introduction to Nonlinerar
Hydroelasticity}, Nauka, 1991, p. 200.


\end{thebibliography}
\end{document}